\documentclass[letterpaper, 10pt, conference]{ieeeconf}
\IEEEoverridecommandlockouts
\usepackage{cite}
\usepackage{amsmath,amssymb,amsfonts,mathtools}
\usepackage{algorithmic}
\usepackage[short]{optidef}
\usepackage{graphicx}
\usepackage{textcomp}
\usepackage{xcolor}
\usepackage{booktabs}
\usepackage{multirow}
\newcommand{\ra}[1]{\renewcommand{\arraystretch}{#1}}

\renewcommand{\arraystretch}{1.2}
\def\BibTeX{{\rm B\kern-.05em{\sc i\kern-.025em b}\kern-.08em
    T\kern-.1667em\lower.7ex\hbox{E}\kern-.125emX}}

\begin{document}

\title{An Implicit and Explicit Dual Model Predictive Control Formulation for a Steel Recycling Process}

\author{{Andrea Ghezzi$^1$, Florian Messerer$^1$, Jacopo Balocco$^2$, Vincenzo Manzoni$^2$, Moritz Diehl$^{1,3}$}
\thanks{$^1$ Department of Microsystems Engineering (IMTEK), University of Freiburg, 79110 Freiburg, Germany \tt{\small\{andrea.ghezzi, florian.messerer, moritz.diehl\}@imtek.uni-freiburg.de}}
\thanks{$^2$ Tenaris Dalmine S.p.A., Data Science Department, Dalmine (BG), Italy}
\thanks{$^3$ Department of Mathematics, University of Freiburg, 79104 Freiburg, Germany}
\thanks{This research was supported by DFG via Research Unit FOR 2401 and project 424107692 and by the EU via ELO-X 953348.
We thank Armin Nurkanović for the helpful comments.
}}
\maketitle

\begin{abstract}
    We present a formulation for both implicit and explicit dual model predictive control for a steel recycling process.
    The process consists in the production of new steel by choosing a combination of several different steel scraps with unknown pollutant content.
    % Here, new steel needs to be produced by choosing a combination of several different steel scraps with unknown pollutant content.
    The pollutant content can only be measured after a scrap combination is molten, allowing for inference on the pollutants in the different scrap heaps.
    The production cost should be minimized while ensuring high quality of the product through constraining the maximum amount of pollutant.
    The dual control formulation allows to achieve the optimal explore-exploit trade-off between uncertainty reduction and cost minimization for the examined problem.
    Specifically, the dual effect is obtained by considering the dependence of the future pollutant uncertainties on the scrap selection in the predictions.
    % By considering the dependence of the future pollutant uncertainties on the scrap selection, the explore-exploit trade-off between uncertainty reduction and cost minimization arises and hence the dual control effect.
    % \todo{Awareness of uncertainty and dual control effect is achieved by including in the predictions the uncertainty on the inference of the pollutants.
    % For numerical stability, these predictions are obtained from a square-root Kalman filter update based on a QR decomposition.
    % In the implicit formulation, the incentive for uncertainty reduction is given indirectly via the impact of active constraints on the objective, as large uncertainty leads to large safety backoffs from the constraint set boundary.
    % The explicit formulation additionally uses a heuristic cost term on uncertainty to encourage its active exploration.
    The implicit formulation promotes uncertainty reduction indirectly via the impact of active constraints on the objective, while the explicit formulation adds a heuristic cost on uncertainty to encourage active exploration.
    We compare the formulations by numerical simulations of a simplified but representative industrial steel recycling process.
    The results demonstrate the superiority of the two dual formulations with respect to a robustified but non-dual formulation.
    % Specifically we achieve a lower cost for the closed-loop trajectories while ensuring constraint satisfaction with a given probability.
% Here, new steel needs to be produced by choosing a combination of several different steel scraps with unknown pollutant content.
% The pollutant content can only be measured after a scrap combination is molten, allowing for inference on the pollutants in the different scrap heaps.
% The cost should be minimized while ensuring high quality of the product through constraining the maximum amount of pollutant.
% The numerical simulations demonstrate the superiority of the two dual formulations with respect to a robustified but non-dual formulation.
% Specifically we achieve lower cost for the closed-loop trajectories while ensuring constraint satisfaction with a given probability.
\end{abstract}

\section{Introduction}
Model-based control approaches such as model predictive control (MPC)~\cite{Rawlings2017} use a model of the controlled system to compute appropriate control inputs.
By operating in a closed-loop framework they are able to react to uncertainties and model-plant mismatch.
However, standard, i.e., nominal MPC has no model of the uncertainty.
Feedback is applied based on the current state estimate independent of its quality. Further, in the absence of additional measures, it will often plan trajectories right on the boundary of the feasible set, such that a slight perturbation can lead to constraint violation.
Approaches such as stochastic or robust MPC (SMPC resp. RMPC) try to solve this problem by explicitly taking into account the uncertainty of the model predictions~\cite{Rawlings2017, Kouvaritakis2016, Mesbah2016, Rakovic2019}.
However, typically this model of uncertainty is static in the sense that they are not aware how uncertainty can be reduced by learning about the system.
The field of dual control~\cite{Feldbaum1960} considers that controls can be used to achieve two competing aims, in what is often referred to as explore-exploit trade-off: exploiting the already available information to (greedily) advance the original control objective, or exploring the system by exciting it in such a way that new information is efficiently obtained.
Model-based approaches can be made aware of this possibility of uncertainty reduction by including an estimator model and planning over the corresponding closed-loop policies~\cite{BarShalom1974}.
As this problem is generally intractable, only approximations can be solved.
If the approximation maintains the relevant dual control aspects it is considered implicit dual control.
Otherwise additional measures need to be taken, such as a heuristic cost term on uncertainty or random control perturbations ensuring sufficient excitation, in what is considered explicit dual control~\cite{Filatov2000}.
An excellent survey on dual MPC can be found in~\cite{Mesbah2018}.

Steel is mainly produced using one of two methods: Blast Furnace (BF) or Electric Arc Furnace (EAF).
The former uses iron ore and cooked coal as raw materials, while the latter melts steel scrap with electrical current.
Steel manufacturing accounts for around 25\% of industrial greenhouse gas emissions~\cite{Allwood2010}, and steelmaking from scrap in an EAF generates one-third of the emissions associated with steelmaking in a BF~\cite{Yellishetty2011}.
% Steel is one of the world most recycled materials, with end-of-life recovery rate estimates as high as 90\%.
Therefore, further reliance on steel scrap recycling and green electric energy is fundamental to achieve a sustainable steel production.
Nonetheless steelmaking via EAF presents challenges.
In fact when steel scraps are molten, we obtain low concentrations of residual elements which are not intentionally added during the steel production cycle and are difficult to remove.
These residual elements may harm the steel properties.
As showed in~\cite{Daehn2017}, steel production from scrap will become even more difficult in the future, as the content of residual elements is likely to increase if less steel is produced from iron ore.
In this work we limit our concern to copper which is a residual element widely spread in mechanical and electrical waste and is known to cause surface defects during hot rolling processes~\cite{Stephenson1983}, thus limiting the applicability of recycled steel.
% one of the main residual elements in mechanical and electrical waste and is known to cause surface defects during hot rolling processes \cite{Stephenson1983}, thus limiting the applicability of recycled steel.
For all these reasons having a method that monitors the residual elements and selects the scrap to melt according to the steel to be produced has paramount importance.

We consider a scrapyard where the scrap is divided into heaps according to their provenience such as automotive or rail industry and their characteristics such as stainless steel, high or low alloy steel scrap.
The aim is to produce steel, which has a maximum content of copper allowed, with the cheapest scrap mix so as to minimize the cost associated with raw materials.
However, given the dimension of the heaps and heterogeneity of the scrap it is not
possible to know exactly the copper content in each heap.
Therefore, we assume to have statistical knowledge of each heap and can refine this knowledge during the production process.
Indeed, once it is decided how much scrap is picked from each heap, the scraps are molten and a measurement
of the copper content in the product can be taken.
The explore-exploit trade-off arises in this scrap selection problem because on the one hand we would like to exploit the current knowledge to achieve the cheapest possible scrap combination, favoring a repetition of previously tried successful selections.
On the other hand changing the scrap mix allows us to improve our knowledge via exploration and possibly leads to a more economical selection in the long run.

\textit{Related work: }
A similar (non-square root based) MPC formulation is proposed in~\cite{Bonzanini2021}, but for a deterministic linear system and with the constraint affected by a state- and control-dependent linear random process.
% \andrea{
%     The dual MPC formulations proposed in~\cite{Heirung2017} and \cite{Soloperto2019} are devoted to output tracking problems for discrete and time-invariant linear systems with output parametrized by an uncertain parameter.}
The dual MPC formulations proposed in~\cite{Heirung2017} and~\cite{Soloperto2019} deal with discrete and time-invariant linear systems with deterministic state but output parametrized by an uncertain parameter.
Although our model can be rewritten in their framework, they consider output tracking problems and explicitly account for the parameter uncertainty in the cost function.
% \andrea{The parameter uncertainty is explicitly accounted in the cost function, the formulation is not square root based and} \notsure{the parameter has no dynamic.}
% for solving an output tracking problem consider discrete and time-invariant linear system with noise-corrupted output but deterministic state.
In Belief-space planning~\cite{Platt2010} an estimator model is used to predict future estimation uncertainty but without considering uncertain constraints.
In~\cite{Falanga2018} the authors introduce a formulation for perception-aware MPC of a quadcopter, but without modelling uncertainty explicitly and instead using cost terms and constraints to keep the nominal state in regions with good observability.
% \andrea{
% In~\cite{Berkenkamp2021} they develop a method for tuning the parameters of a generic nonlinear system to maximize a performance measure subject to safety constraints where the dependance of the parameters on the performance measure is not known a priori and the safety constraints are uncertain. Their algorithm start from a given safe parametrization and try to find the largest safe set (exploration phase) while computing the optimal parameter for the current safe (exploitation phase). Gaussian Processes are used to propagate uncertainty along the system trajectory, which make difficult to apply the algorithm on high-dimensional or embedded system due to the demanding computational power.
% }
% \notsure{I might remove the reference to Berkenkamp or write it differently - difficult to shorten the approach to a sentence since it's quite different compared to all the others.\\}
% \notsure{Reference to bandit problem with stochastic constraints?}

% \paragraph{Contribution}
\textit{Contribution: }
In this paper we introduce an implicit and explicit dual MPC formulation for a steel recycling process.
The advantages of the proposed formulations against non-dual approaches are demonstrated in numerical experiments.
The simulation is based on fictitious numbers and reduced in dimension for clarity of exposition, but otherwise realistic to how it could be used on a real plant.
The process is modelled as an autonomous linear system where the state is estimated with a Kalman filter.
System control and dual effect occur when the augmented state formed by the state estimate and its covariance is predicted, since the covariance is directly affected by the scrap selection.
A QR-factorization based square-root Kalman filter update improves the numerical stability and ensures positive-semidefinitness of the predicted covariance matrices as compared to a standard Kalman filter.

\textit{Outline: }
In Sec.~\ref{sec: problem formulation} we introduce the model of the problem together with the nominal and robustified formulation.
In Sec.~\ref{sec: dual formulation} we present the dual formulations and, eventually, in Sec.~\ref{sec: numerical examples} we compare the formulations by numerical simulations of a simplified steel recycling process.

\section{Nominal and robust formulation}\label{sec: problem formulation}
% \subsection{System Modelling}
% We consider a discrete-time uncertain system that is affine in the uncertain state and disturbances but possibly nonlinear in the control input:
% \begin{equation} \label{eq: generical non linear system}
%     \begin{cases}
%         x_{t+1} &= A(u_t) x_t + B(u_t) w_t + r(u_t)\\
%         y_t &= C(u_t) x_t + D(u_t) v_t + s(u_t)
%     \end{cases}, \quad t=0, \dots, T-1,
% \end{equation}
% where $x_t \in \mathbb{R}^{n_x}$ denotes the state at discrete time, $y_t \in \mathbb{R}^{n_y}$ the output and they can be influenced by unknown i.i.d. disturbances $w_t \in \mathbb{R}^{n_w}$ and $v_t \in \mathbb{R}^{n_v}$ respectively.
% Note that $A: \mathbb{R}^{n_u} \rightarrow \mathbb{R}^{n_x \times  n_x}$, $B: \mathbb{R}^{n_u} \rightarrow \mathbb{R}^{n_x \times  n_w}$, $C: \mathbb{R}^{n_u} \rightarrow \mathbb{R}^{n_y \times n_x}$, $D: \mathbb{R}^{n_u} \rightarrow \mathbb{R}^{n_y \times n_v}$, $r: \mathbb{R}^{n_u} \rightarrow \mathbb{R}^{n_x}$ and $s: \mathbb{R}^{n_u} \rightarrow \mathbb{R}^{n_y}$ are general nonlinear functions in $u$ that we assume differentiable.
We model the steelmaking process as an autonomous discrete time linear system where time index $t=0,\dots,T$ denotes the casts' sequence.
The model is given by
\begin{equation} \label{eq: linear system}
    \begin{cases}
        x_{t+1} &= x_t + w_t \\
        y_t &= u_t^\top x_t + v_t
    \end{cases}, \quad t=0, \dots, T,
\end{equation}
where the state $x_t \in \mathbb{R}^{n_x}$ represents the copper content in each heap of scrap, the output $y_t \in \mathbb{R}^{n_y}$ contains the measured copper content in the produced steel.
The output is a linear combination of the state weighted by the controlled variable $u_t \in \mathbb{R}^{n_x}$ which corresponds to the amount of scrap picked from each heap.
State and output are corrupted by normally distributed and i.i.d.~disturbances $w_t \sim \mathcal{N}(0, Q)$ and $v_t \sim \mathcal{N}(0, R)$ respectively, where $\mathcal{N}(\mu, \Sigma)$ denotes the normal distribution with mean $\mu$ and covariance $\Sigma$.

% Regarding the steelmaking process, the time index $t=0,\dots,T-1$ denotes the casts' sequence, the state represents the copper content in each heap of scrap, the control represents the amount of scrap picked from each heap and the output the measured copper content in the produced steel.
% We model the system dynamics as a random walk process.
% Thus, in our case $A(u_t) \coloneqq I$, $B(u_t) \coloneqq \mathbf{1}$, $w_t \sim \mathcal{N}(0, Q)$ and i.i.d.,
% $r(u_t) \coloneqq \mathbf{0}$, $s(u_t) \coloneqq 0$ where $I$, $\mathbf{1}$ and $\mathbf{0}$ denote respectively the identity matrix, the unitary matrix and the zero vector with corresponding dimensions, and $\mathcal{N}(\mu, \Sigma)$ denotes the Gaussian distribution with mean $\mu$ and covariance $\Sigma$.
% Finally, the resulting model is an autonomous discrete-time linear system given by:
% where the initial state $x_0$ is assumed to be normally distributed $\mathcal{N}(\hat{x}_0, P_0),$ with covariance $P_0 \succ 0$.

In the considered set-up the system state $x_t$ is not accessible, therefore a state observer is needed to compute a state estimate $\hat{x}_t$, we assume that the initial state $x_0$ is normally distributed $\mathcal{N}(\hat{x}_0, P_0)$ with covariance $P_0 \succ 0$.
Despite the linear system in itself is autonomous, when we consider the estimated state $\hat{x}_t$ augmented by its covariance $P_t$, the control variable $u_t$ directly affects $P_t$ via the estimator update, therefore obtaining a controlled system.
% We then obtain a peculiar form of control where we achieve the control of the system output $y_t$ by manipulating $u_t$.
% In fact, when we consider the estimated state $\hat{x}_t$ augmented by its covariance $P_t$, the control variable $u_t$ directly affects $P_t$ via the estimator update.
% Thus, even though the linear system in itself is autonomous, the estimation covariance is affected by the controls.

\subsection{Optimal scrap selection problem -- Nominal formulation}
We assume that every scrap has a certain price $p_i, i=1, \dots, n_x$ and the concentration of copper, i.e., the pollutant, in the final steel $y_t$ is upper bounded by a value $y_{\max}$. Therefore we aim to minimize the cost of the steel produced in the cast $t$ by selecting the cheapest mix of scraps $u_t$ which fulfills the limitation on $y_t$.

A first simple formulation of the scrap selection problem is given by \\
\parbox{0.18\columnwidth}{
   \raggedleft
   $\mathcal{P}_\mathrm{NOM}:$}
\hfill
\parbox{.8\columnwidth}{
\begin{mini!}[2]
    {u_0}{p^\top u_0}
    {\label{ocp: naive formulation}}{}
    \addConstraint{u_0^\top \hat{x}_0}{\leq y_{\max}\label{cns: robust cns naive formulation}}
    \addConstraint{u_{\min} \leq u_0}{\leq u_{\max}}
    \addConstraint{\mathbf{1}^\top u_0}{= 1,}
\end{mini!}
}
where $\hat{x}_0$ is the state estimate at the current time instant $t$ obtained from a Kalman filter.
Note that this formulation is not accounting for uncertainty on the state.
Thus, if $\hat{x}_0$ is a wrong guess of $x_0$, constraint~\eqref{cns: robust cns naive formulation} might be violated by the real state $x_0$.
As in this formulation the state is not affected by the controls, there is no reason to have a prediction horizon.
Finally, note that the dynamics~\eqref{eq: linear system} plays a role only in closed-loop control.
Indeed, we apply the optimal $u_0$ to the process and measure the corresponding $y_0$.
Then, we feed these two quantities to a Kalman filter and update our estimate of $\hat{x}_0$.
We will discuss the closed-loop scheme at the end of Section~\ref{sec: dual formulation}.
% Moreover, problem~\eqref{ocp: naive formulation} doesn't entail a prediction horizon.

% Since the state $x$ is not accessible we can introduce an estimation routine, e.g., based on a Kalman filter, which gives us
% an estimate of the state with mean $\hat{x}$ and covariance $P$.
\subsection{Robust formulation}
Considering the covariance $P$ of the state estimate $\hat{x}$, we can take into account its level of uncertainty.
Thus, we can require constraint~\eqref{cns: robust cns naive formulation}
to hold with at least a certain probability, resulting in the chance constraint
\begin{equation}\label{eq: chance constraint}
    \mathbb{P}_{x \sim \mathcal{N}(\hat{x}, P)}(u^\top x \leq y_{\max}) \geq 1 - \epsilon,
\end{equation}
with $\epsilon>0$ being the maximal allowed probability of constraint violation.
Since $x$ is normally distributed, the non-noisy output $u^\top x$ follows a normal distribution as well with $u^\top x \sim \mathcal{N}(u^\top \hat{x}, u^\top P u)$.
Therefore, we can write the chance constraint as
\begin{equation}\label{eq: explicit version of the chance constraint}
    u^\top \hat{x} + \gamma(\epsilon) \sqrt{u^\top P u} \leq y_{\max},
\end{equation}
where the coefficient $\gamma(\epsilon) = \Phi^{-1}(1-\epsilon)$ controls the safety backoff from the constraint and can be computed from the inverse cumulative distribution function $\Phi^{-1}$ of the standard normal distribution
such that the amount of backoff corresponds to the specified probability $1-\epsilon$.
In the following, we drop the explicit dependency of $\gamma$ on $\epsilon$.
We also refer to this constraint as robustified because for a fixed value of $\gamma$ it corresponds to a robust constraint where a bounded distribution of the state is assumed, with support given by the ellipsoidal level line of the normal distribution corresponding to the chosen value of $\gamma$.
The resulting robustified scrap selection problem can be formulated as a Second Order Cone Program (SOCP)~\cite{Boyd2004}: \\
\parbox{0.18\columnwidth}{
   \raggedleft
   $\mathcal{P}_\mathrm{ROB}:$}
\hfill
\parbox{.8\columnwidth}{
\begin{mini!}[2]
    {u_0}{ p^\top u_0}
    {\label{ocp: stochastic formulation}}{}
    \addConstraint{u_0^\top \hat{x}_0 + \gamma \sqrt{u_0^\top P_0 u_0}}{\leq y_{\max}\label{cns: robust cns stochastic formulation}}
    \addConstraint{u_{\min} \leq u_0}{\leq u_{\max}}
    \addConstraint{\mathbf{1}^\top u_0}{= 1.}
\end{mini!}
}
As in the nominal case, there is no prediction horizon, since the influence of the scrap selection on future uncertainty is not modelled.

\section{Dual formulations}\label{sec: dual formulation}
In this section we illustrate the dependency of the state estimate covariance on the scrap selection and how the dual effect is obtained from the covariance predictions.
% In the next section we show how to introduce this relation in the optimal scrap selection problem and obtain the dual effect.
%  the this state estimation covariance can provide a feedback for the controlled variable $u$.
\subsection{Implicit dual formulation}
Consider the Kalman filter propagation of the state estimate covariance, the future covariance is directly affected by $u_t$, then we can manipulate future uncertainty.
% However, if we consider the Kalman filter propagation of the covariance matrix, the future covariance is affected
% by $u_t$, then we can manipulate future uncertainty.
Hence, it is meaningful to consider a prediction horizon $k=0,\dots,N$ for the optimal scrap selection problem where the propagation of the state estimate covariance anticipate feedback for $u$.
Since the chance constraint needs to hold when taking the measurement, i.e., without taking into account the new information, the predicted covariance $P_k \coloneqq P_{k|k-1}$ is the relevant one.
The resulting optimization problem is given by
\begin{mini!}[3]
    {\substack{u_0, \dots, u_N, \\ P_{1}, \dots, P_{N}, \\
               K_0, \dots, K_{N-1}}}{\sum_{k=0}^{N} p^\top u_k}
    {\label{ocp: multiple shooting formulation}}{}
    \addConstraint{P_{k+1}}{= \psi(u_{k}, K_{k}, P_{k}), \;}{k=0,\dots,N-1}
    \addConstraint{K_k}{= P_{k} u_k (u_k^\top P_{k} u_k + R)^{-1}, \;}{ k=0,\dots,N-1}
    \addConstraint{u_k^\top \hat{x}_0 + \gamma \sqrt{u_k^\top {P}_{k} u_k}}
                  {\leq y_{\max}, \;}{k=0,\dots,N\label{cns: robust cns dual ms formulation}}
    \addConstraint{u_{\min} \leq u_k }{\leq u_{\max}, \;}{k=0,\dots,N}
    \addConstraint{\mathbf{1}^\top u_k}{= 1, \;}{k=0,\dots,N,}
\end{mini!}
where $K_k$ are the Kalman gains and function $\psi$ denotes the covariance propagation
\begin{equation}\label{eq: joseph form cov prop}
    \begin{array}{l}
        \psi(u_k, K_k, P_k) = \\
        = (I - K_k u_k^\top) P_k (I - K_k u_k^\top)^\top + K_k R K_k^{\top} + Q.
    \end{array}
\end{equation}
Recall that $Q$ and $R$ are the covariance of the state and output noise respectively.
Note that the chance-constraint~\eqref{cns: robust cns dual ms formulation} is enforced with respect to $\hat{x}_0$ since due to the model the expectation of the state remains constant throughout the horizon.
Since the propagation of $P_{k+1}$ depends on $u_{k}$, the optimization problem~\eqref{ocp: multiple shooting formulation} can act on the backoff term in constraint~\eqref{cns: robust cns dual ms formulation}.
Therefore, the incentive to reduce the uncertainty depends exclusively on the reduction of the backoff and the exploration incentive increases with the horizon length as a longer horizon allows for more exploitation.
Intuitively, if the horizon has length $N=1$, there is no chance to explore because the optimization problem can only exploit the available information for selecting the scrap.
On the other hand, if the prediction horizon is sufficiently long, it can be worthwhile to take non-greedy actions that reduce the future uncertainty.
Thus, the formulation~\eqref{ocp: multiple shooting formulation} is not only statically uncertainty-aware but it is also able to deliberately reduce future uncertainty.
It is a strongly nonlinear and nonconvex optimization problem that can be solved with a NLP solver.

Even though in~\eqref{eq: joseph form cov prop} we update the covariance matrix according to the \textit{Joseph form}, which retains positive definiteness and symmetry of the covariance matrix~\cite{Stengel1994}, this is not guaranteed during the optimization solver iterations, resulting in possible numerical difficulties.
It is possible to overcome this issue by considering the propagation of the square root factors of the covariance matrix, as shown in the following.

\subsection{Square root covariance propagation}
Following the square root covariance filtering algorithm introduced in~\cite{Morf1975} it is possible to propagate the
covariance matrix just using the corresponding square root matrices.
Given a positive definite matrix $A \succ 0$, a square-root factor will be defined as any matrix, $A^\mathrm{r}$, such that $A = {A^{\mathrm{r}}}^\top A^\mathrm{r}$.
In general, square root factors are not unique.
They can be made unique by imposing specific properties such as symmetry or triangular structure. In our case the latter is preferred, therefore in the following any
$A^\mathrm{r}$ denotes the upper triangular factor.

Denote the upper triangular decomposition of the matrices $P_0, Q, R$ as $P_0^{\mathrm{r}}, Q^\mathrm{r}, R^\mathrm{r}$.
Then for $k=0, \dots, N-1$ we can propagate $P_k^{\mathrm{r}}$ by a QR decomposition~\cite{Morf1975} as follows
\begin{equation} \label{eq: sqrt cov update}
    \begin{bmatrix}
        R^\mathrm{r} & 0 \\
        \boxed{P_{k}^{\mathrm{r}} u_k} & \boxed{P_{k}^\mathrm{r}}\\
        0 & Q^\mathrm{r}
    \end{bmatrix} =
    Q_k^\mathrm{F}
    \begin{bmatrix}
        S_k^\mathrm{r} & P_{k} u_k S_k^{- \mathrm{r}} \\
        0 & \boxed{P_{k+1}^\mathrm{r}}\\
        0 & 0
    \end{bmatrix},
\end{equation}
where $Q_k^\mathrm{F}$ is an orthogonal matrix and $S_k^\mathrm{r}$ denotes the square root factor of the observer innovation covariance $S_k = u_k^\top P_k u_k + R = S_k^\mathrm{r \top} S_k^\mathrm{r}$.
Then,~\eqref{eq: sqrt cov update} defines a function $\psi_{\mathrm{QR}}$ such that 
\begin{equation}
    P_{k+1}^\mathrm{r} = \psi_{\mathrm{QR}}(u_k, P^\mathrm{r}_k).
\end{equation}
Moreover, it is possible to operate only with the square root factors of the covariance matrix to preserve their condition number throughout the prediction horizon of problem~\eqref{ocp: multiple shooting formulation}.
Therefore, the term $\sqrt{u_k^\top P_{k} u_k}$ in constraint~\eqref{cns: robust cns dual ms formulation} is restated as
\begin{equation}
    \sqrt{(P_{k}^\mathrm{r} u_k)^\top(P_{k}^\mathrm{r} u_k)} =
    {\lVert P_{k}^\mathrm{r} u_k  \rVert}_2.
\end{equation}

Finally, the optimal scrap selection problem~\eqref{ocp: multiple shooting formulation} can be equivalently stated as $\mathcal{P}_\mathrm{IMPL}$:
\begin{mini!}[3]
    {u_0, \dots, u_N}{\sum_{k=0}^{N} p^\top u_k}
    {\label{ocp: single shooting sqrt formulation}}{\label{cf: single shooting sqrt formulation}}
    \addConstraint{u_k^\top \hat{x}_0 + \gamma {\lVert \tilde{P}_{k}^\mathrm{r} (u) u_k \rVert}_2}
                    {\leq y_{\max}, \quad}{k=0,\dots,N\label{cns: robust cns dual ss formulation}}
    \addConstraint{u_{\min} \leq u_k }{\leq u_{\max}, \quad}{k=0,\dots,N}
    \addConstraint{\mathbf{1}^\top u_k}{= 1, \quad}{k=0,\dots,N,}
\end{mini!}
where $\tilde{P}_{0}^\mathrm{r}$ corresponds to the Cholesky decomposition of $P_0$, $u=(u_0, \dots, u_N)$ and $\tilde{P}_{k+1}^\mathrm{r}(u) = \psi_{\mathrm{QR}}(\tilde{P}^{\mathrm{r}}_{k}(u), u_{k}), k=0,\dots, N-1$.
In contrast to the formulation~\eqref{ocp: multiple shooting formulation}, here the covariance propagation is not expressed as constraints but is computed externally as function of the controls $u$.

\subsection{Explicit dual formulation}
By adding a heuristic cost on uncertainty it is possible to encourage exploration even more than via the implicit incentive provided by the backoffs in~\eqref{cns: robust cns dual ss formulation}.
This new term is weighted in the cost function by a hyperparameter $\alpha$ that regulates the explore-exploit trade-off.
Therefore, the explicit dual formulation $\mathcal{P}_\mathrm{EXPL}$ share the same constraints of formulation~\eqref{ocp: single shooting sqrt formulation} but its cost function is given by
\begin{equation}\label{eq: explicit exploration cost function}
    \sum_{k=0}^{N} p^\top u_k + \alpha \mathrm{Tr}(\tilde{P}_{k}^\mathrm{r}(u)^\top \tilde{P}_{k}^\mathrm{r}(u)),
\end{equation}
where $\alpha \geq 0$ is the hyperparameter that tunes the exploration incentive.
Adding the trace of the covariance matrix to the cost function allows us to evenly minimize uncertainty along all directions.
Note that for $\alpha=0$ we recover the cost function~\eqref{cf: single shooting sqrt formulation} of the implicit dual MPC problem.

\subsection{Closed-loop scheme}
\begin{figure}
    \vspace{0.3cm}
    \centering
    \includegraphics[width=\linewidth]{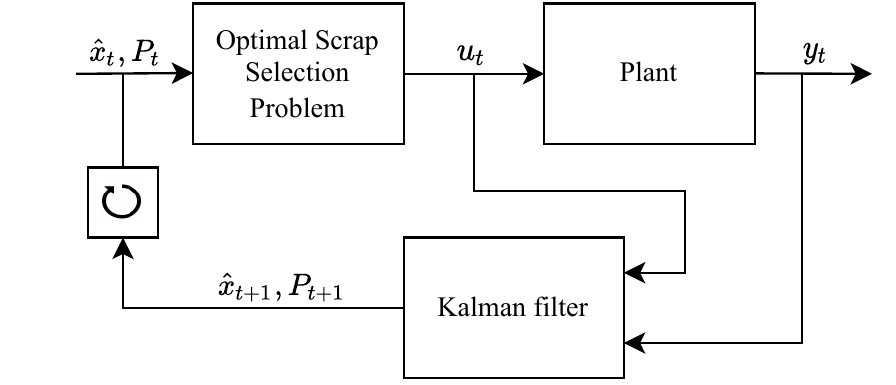}
    \caption{Closed loop block scheme.}\label{fig: closed-loop scheme}
\end{figure}

We are interested in solving the optimal scrap selection problems~\eqref{ocp: naive formulation},~\eqref{ocp: stochastic formulation},~\eqref{ocp: single shooting sqrt formulation} and~\eqref{ocp: single shooting sqrt formulation} with cost function~\eqref{eq: explicit exploration cost function} for every cast $t=0,\dots,T-1$.
We apply to the plant only the first element of solution vector $u_t \coloneqq u^\star_0$ following the
receding horizon principle.
Once the cast $t$ is completed we can measure the copper concentration $y_t$ and feed both $u_t$ and $y_t$ to the Kalman filter routine.
In this way, we update the latest predictions of $\hat{x}_{t+1}, P_{t+1}$, needed to solve the optimization problem at the next time step $t+1$.
A block scheme for the closed-loop simulation is depicted in Fig.~\ref{fig: closed-loop scheme}.

\section{Numerical examples}\label{sec: numerical examples}
In this section, we compare the closed-loop trajectories of the problem formulations stated in Sec.~\ref{sec: problem formulation}-\ref{sec: dual formulation}.
First, we show an example for each formulation and then, since we have an uncertain system, we assess the general closed-loop behavior by sampling many different uncertainty realizations.

\begin{table}
    \centering
    \ra{1.2}
    \caption{Simulation parameters used in Sec.~\ref{sec: numerical examples}}
    \begin{tabular}{@{}lrr@{}} \toprule
        Name & Symbol & Value \\ \midrule
        True initial state & $x_0$ & $(0.07, 0.13, 0.17)$ \\
        Initial state covariance & $P_0$ & $\mathrm{diag}(10^{-4}, 10^{-3}, 10^{-3})$ \\
        State noise covariance & $Q$ & $10^{-7} \cdot I$ \\
        Output noise covariance & $R$ & $2\cdot10^{-6}$ \\
        Scrap prices & $p$ & $(2,1,1)$ \\
        Max copper allowed & $y_{\max}$ & $0.12$ \\
        Max constraint violation & $\epsilon$ & $2.55\%$ \\
        Max constraint violation & $\gamma(\epsilon)$ & $2$ \\
        Control constraints & $[u_{\min}, u_{\max}]$ & $[0, 1]$ \\
        Exploration hyperparameter & $\alpha$ & $100$ \\
        Prediction horizon & $N$ & $15$ \\
        Simulation length & $T$ & $20$ \\ \bottomrule
    \end{tabular}\label{tab: parameters of selected example}
\end{table}

The four formulations share the same parameters and initialization, which are collected in Table~\ref{tab: parameters of selected example}.
% The four formulations share the same parameters and initialization. Specifically, the state dimension is $n_x = 3$.
% we suppose to know the initial state $x_0 = (0.07, 0.13, 0.17)$, while the initial estimated state $\hat{x}$ is sampled 
% from $\mathcal{N}(x_0, P_0)$ where $P_0 = \mathrm{diag}(10^{-4}, 10^{-3}, 10^{-3})$.
The true initial state of the system $x_0$ is always the same, but we assume imperfect knowledge of it. 
In consequence the initial state estimate $\hat x_0$ varies, e.g., because of a different history, and is sampled as $\hat x_0 \sim \mathcal{N}(x_0, P_0)$.
% The state disturbance $w_t$ has covariance $Q=10^{-7} \cdot I$ and the output noise $v_t$ has variance $R=2\cdot10^{-6}$.
% The prices of the three scraps are given by $p=(2,1,1)$.
The choice of the scrap prices $p$ is motivated by the fact that the first scrap has higher cost since its copper content is lower in terms of mean value and uncertainty compared to the other two.
% We aim to produce a steel with a maximum allowed copper content $y_{\max} = 0.12 \%$ and in the uncertainty-aware formulations we set the backoff coefficient to $\gamma=2$, i.e., we allow a maximal probability violation of the chance-constraint~\eqref{eq: chance constraint}  $\epsilon \simeq 2.55\%$.
% Regarding the scrap availability, we suppose that each heap can supply an infinite amount of scrap and in each formulation we impose $u_{\min} =0$ and $u_{\max}=1$.
Moreover, we suppose that each heap can supply an infinite amount of scrap, so that we will not run out of a scrap during the closed-loop simulation.
The choice of the exploration hyperparameter $\alpha$ will be discussed at the end of the Section.
% This choice will be motivated at the end of the Section.
% All the closed-loop simulations have length $T=20$.
The simulations are carried out using Python, the optimization problems are formulated using CasADi~\cite{Andersson2019} and solved via IPOPT~\cite{Waechter2006}.

In the following we refer to the formulations~\eqref{ocp: naive formulation},~\eqref{ocp: stochastic formulation},~\eqref{ocp: single shooting sqrt formulation} and~\eqref{ocp: single shooting sqrt formulation} with cost function~\eqref{eq: explicit exploration cost function} as \textit{nominal}, \textit{robust}, \textit{implicit dual} and \textit{explicit dual} formulation, respectively.
% We refer to formulation \eqref{ocp: stochastic formulation} as \textit{robust} formulation because once $\gamma$ is set
% in the chance-constraint \eqref{eq: explicit version of the chance constraint},
% this corresponds to a robust constraint where a bounded distribution of the state is assumed, with support given by the ellipsoidal level line of the normal distribution corresponding to the chosen value of $\gamma$.

% Maybe a TODO:
% \textcolor{blue}{maybe call it ``robustified'' instead.}

The letters $\mathrm{A}, \mathrm{B}, \mathrm{C}$ encode the name of the scrap, therefore $x = (x_\mathrm{A}, x_\mathrm{B}, x_\mathrm{C})$, $\hat{x} = (\hat{x}_\mathrm{A}, \hat{x}_\mathrm{B}, \hat{x}_\mathrm{C})$ and $u = (u_\mathrm{A}, u_\mathrm{B}, u_\mathrm{C})$.

\subsection{Selected example}\label{subsec: numerical examples - selected example}

\begin{figure}
    \vspace*{3mm}
    \begin{center}
        \includegraphics[width=\linewidth]{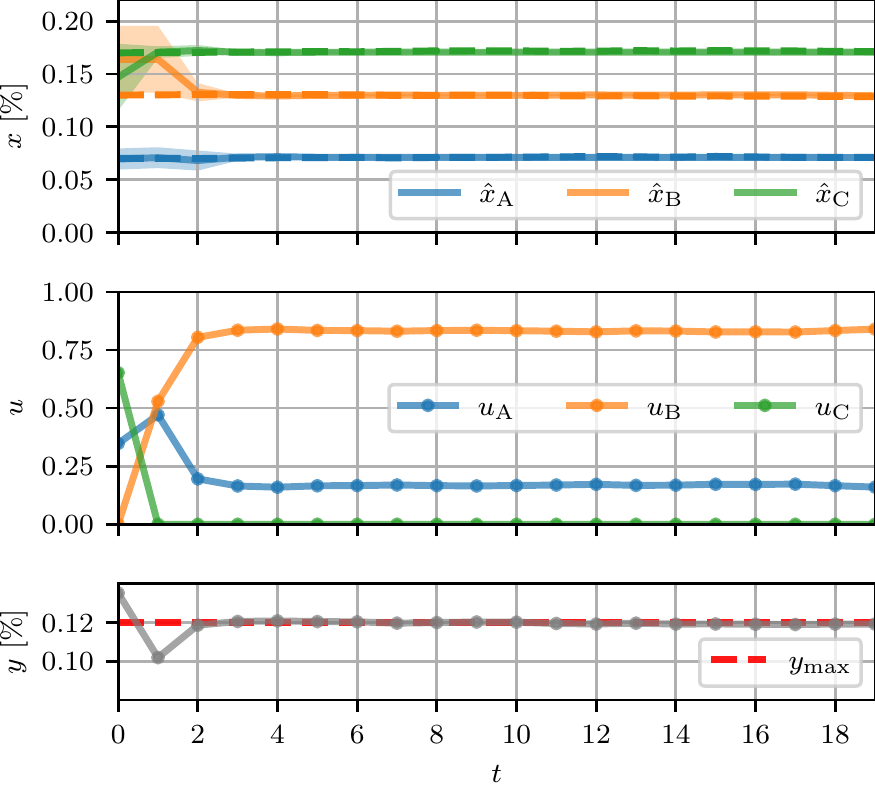}
    \end{center}
    \vspace*{-3mm}
    \caption{Selected example -- \textit{nominal} formulation}\label{fig: selected example - nominal formulation}
\end{figure}

\begin{figure}
    \vspace*{3mm}
    \begin{center}
        \includegraphics[width=\linewidth]{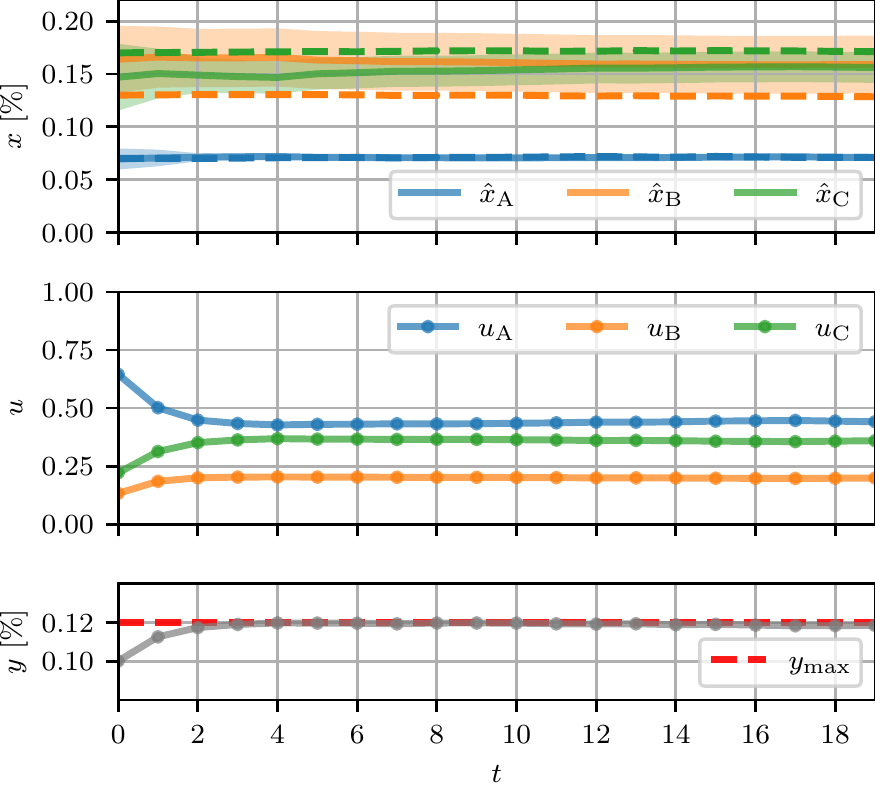}
    \end{center}
    \vspace*{-3mm}
    \caption{Selected example -- \textit{robust} formulation}\label{fig: selected example - stochastic formulation}
\end{figure}

\begin{figure}
    \vspace*{3mm}
    \begin{center}
        \includegraphics[width=\linewidth]{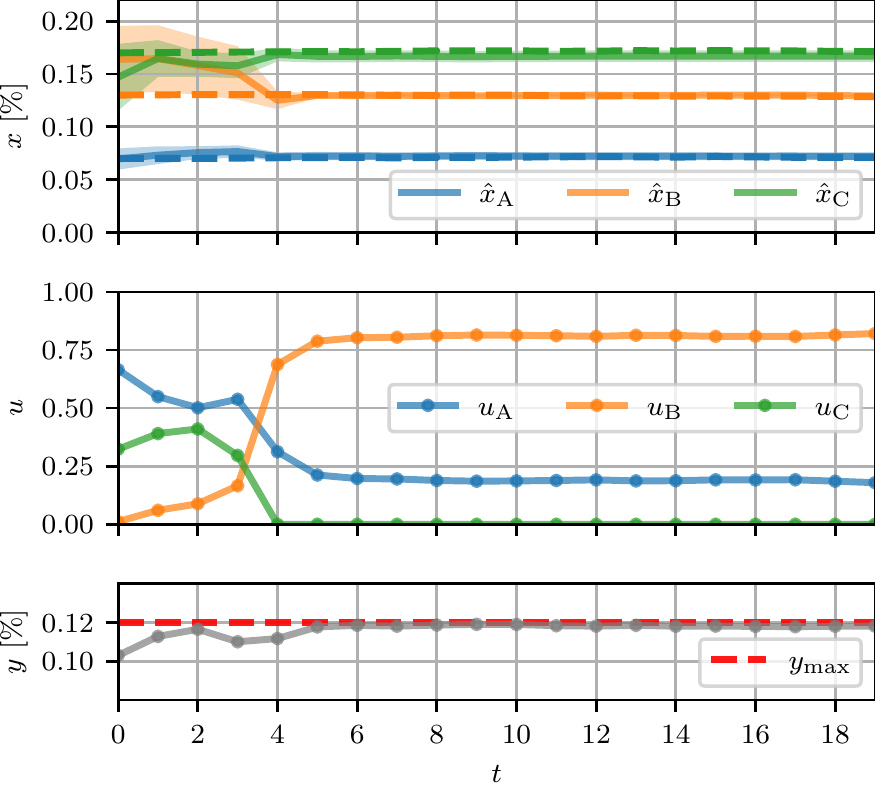}
    \end{center}
    \vspace*{-3mm}
    \caption{Selected example -- \textit{implicit dual} formulation}\label{fig: selected example - implicit dual formulation}
\end{figure}

\begin{figure}
    \vspace*{3mm}
    \begin{center}
        \includegraphics[width=\linewidth]{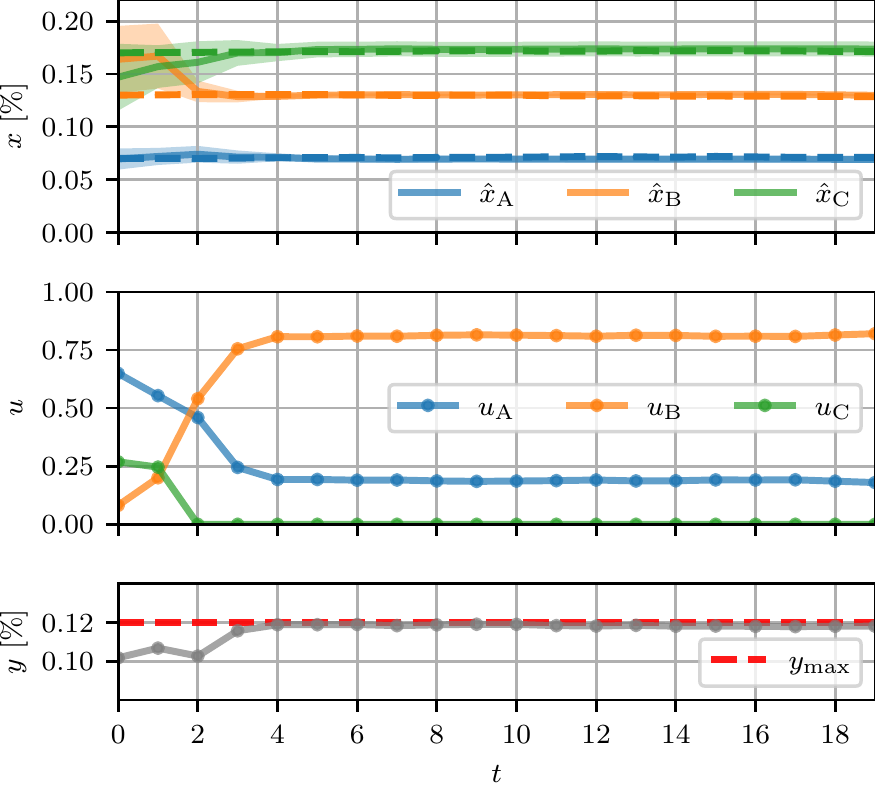}
    \end{center}
    \vspace*{-3mm}
    \caption{Selected example - \textit{explicit dual} formulation}\label{fig: selected example - explicit dual formulation}
\end{figure}

The four formulations share the initial estimated state $\hat{x}_0 = (0.0695 , 0.1639, 0.1469)$ and the same realization of the two disturbances $w_t$, $v_t$, $t=0, \dots, T-1$.
Figures from~\ref{fig: selected example - nominal formulation} to~\ref{fig: selected example - explicit dual formulation} contain the results of the closed-loop simulation.
Specifically, the first one depicts the state, i.e., the copper content of each heap as percentage by weight, the dashed lines are the true state while the solid lines are the estimated state. The shaded area corresponds to the $1\sigma$ standard deviation of $\hat{x}$.
% We carry out the closed-loop simulation according to the parameters described above and for each formulation we show three plots.
% The first one depicts the state, i.e., the copper content of each heap as percentage by weight, the dashed lines are the true state while the solid lines are the estimated state. The shaded area corresponds to the $1\sigma$ standard deviation of $\hat{x}$.
The second plot represents the control, i.e., the selected mass fraction of each scrap.
The third plot pictures the constraint on the maximum copper allowed in the final steel.
% The red dashed line represents $y_{\max}$, while the gray line corresponds to the true or uncorrupted copper content.
This plot shows the possible constraint violations and, for the uncertainty-aware formulations, the magnitude of the backoff.

Fig.~\ref{fig: selected example - nominal formulation} shows the resulting closed-loop trajectory obtained with the \textit{nominal} formulation. One can notice that at $t=1$, the state estimate $\hat{x}_\mathrm{C}$ has recovered the true value $x_\mathrm{C}$.
At this step, $\hat{x}_\mathrm{C}$ is greater than $\hat{x}_\mathrm{B}$, thus the scrap selection picks only scrap $\mathrm{A}$ and $\mathrm{B}$.
This scrap combination strongly excites the Kalman filter block in closed loop, leading to very accurate prediction of the states at time $t=2$.
The accurate estimate of $\hat{x}_\mathrm{B}$ leads to increased use of scrap $\mathrm{B}$ and reduced use of scrap $\mathrm{A}$ to achieve the minimal cost.
This formulation achieves the lowest cost among the ones presented, namely 23.86, but ignores constraint violations, as can be seen in the third plot of Fig.~\ref{fig: selected example - nominal formulation} where the constraint on the maximum copper content allowed is exceeded.
% As shown in Table~\ref{tab: selected example - cost of the closed-loop trajectories} this formulation is the best in terms of cost but ignores constraint violations, as can be seen in the third plot of Fig.~\ref{fig: selected example - nominal formulation} where the constraint on the maximum copper content allowed is exceeded.

Fig.~\ref{fig: selected example - stochastic formulation} is obtained adopting the \textit{robust} formulation to solve the scrap selection problem.
Despite being an uncertainty-aware formulation, it does not know how to actively reduce uncertainty, since there is no explicit dependence of the scrap selection on the state estimate covariance.
From $t=2$ the scrap mix is kept constant for the rest of the simulation. This does not bring new information to the Kalman filter block, resulting in a state estimate $\hat{x}$ far from the true value $x$ and large uncertainty on $\hat{x}_\mathrm{B}$ and $\hat{x}_\mathrm{C}$.
In the end, this leads to a greater cost of the closed-loop trajectory. Yet, from the third plot, one can notice that a backoff is always kept from $y_{\max}$ which avoids constraint violations.

Fig.~\ref{fig: selected example - implicit dual formulation} depicts the closed-loop trajectory attained with the \textit{implicit dual} formulation.
One can see that until $t=4$, scrap $\mathrm{B}$ is barely selected.
Also, $u_\mathrm{C}$ is greater than $u_\mathrm{B}$ until the mean of the estimate $\hat{x}_\mathrm{C}$ is smaller than $\hat{x}_\mathrm{B}$.
However, a small $u_\mathrm{B}$ is picked which is enough to excite the Kalman filter in closed loop such that at $t=5$ we have the opposite situation, $\hat{x}_\mathrm{B}$ smaller than $\hat{x}_\mathrm{C}$.
At this time step, a greater $u_\mathrm{B}$ is selected because this formulation embeds the notion that the scrap selection can reduce the uncertainty on the future estimates leading to lower cost.
The greater use of $u_\mathrm{B}$ has further distinguished $\hat{x}_\mathrm{B}$ and $\hat{x}_\mathrm{C}$ and reduced the uncertainty on the former.
The scrap selection at steps $t=6,7,8$ further increases the mass $u_\mathrm{B}$ while reducing $u_\mathrm{A}$ and eliminating $u_\mathrm{C}$.
From $t=8$ the scrap mix is the same until the end of the simulation. Eventually, the exploration embedded in the problem formulation improves the closed-loop performance compared to the \textit{robust} formulation. This form of exploration is cautious since the realized $y$ is always below the prescribed limit $y_{\max}$, and it involves exclusively the directions where reducing uncertainty leads to lower cost.
The first plot of Fig.~\ref{fig: selected example - implicit dual formulation} shows the main drawback of this approach.
In fact exploration is triggered only when $\hat{x}_\mathrm{B}$ gets smaller than $\hat{x}_\mathrm{C}$ and dependent on the Kalman filter block in closed loop.
This happens because the formulation only plans an open-loop control trajectory instead of a policy.
In consequence, it does not know that the scrap selection will be adapted after new knowledge is acquired.
Instead, with respect to the cost, it only plans with the currently estimated mean.
Note that if the initial state estimates are further apart, we may need too many casts to meet the condition that triggers exploration, leading to high-cost trajectories similar to the ones obtained with the \textit{robust} formulation.

Finally, Fig.~\ref{fig: selected example - explicit dual formulation} shows the results achieved by the \textit{explicit dual} formulation, where exploration is promoted by the additional minimization term in the cost function~\eqref{eq: explicit exploration cost function}.
This formulation reduces the uncertainty in every direction earlier than the previous approach, because it provides scrap combinations which excite more the Kalman filter.
Indeed, already in the scrap mix of step $t=0$, this formulation chooses a bigger $u_\mathrm{B}$ than the \textit{implicit dual} formulation.
Moreover, once uncertainties on $\hat{x}$ are reduced, the second term of cost function~\eqref{eq: explicit exploration cost function} becomes negligible compared to the first one.
Thus, the \textit{explicit dual} formulation doesn't perform any other exploration actions which may increase the cost of the closed-loop trajectories.
This leads to the minimum cost among the uncertainty aware formulations without any constraint violation with a value of 24.92, followed by the implicit dual formulation with 25.42 and the robust one with 29.03.
% As shown in Table~\ref{tab: selected example - cost of the closed-loop trajectories}, this leads to the minimum cost among the uncertainty aware formulations without any constraint violation, as shown in the third plot of Fig.~\ref{fig: selected example - explicit dual formulation}.

% \begin{table}
%     \centering
%     \ra{1.2}
%     \caption{Selected example \\ cost of the closed-loop trajectories}
%     \begin{tabular}{@{}lrrrr@{}} \toprule
%         & \textbf{\textit{Nominal}}& \textbf{\textit{Robust}}& \textbf{\textit{Impl. dual}}& \textbf{\textit{Expl. dual}} \\ \midrule
%         Cost & 23.86 & 29.03 & 25.42 & 24.92 \\ \bottomrule
%     \end{tabular}
%     % \vspace*{-4mm}
%     \label{tab: selected example - cost of the closed-loop trajectories}
% \end{table}

\subsection{Extensive comparison}
Since the performance of each formulation depends on the disturbance and initial state realizations, we run $N_\mathrm{s}=1000$ different simulations, where the adopted parameters are the same as described at the beginning of this section.
The closed-loop trajectories, obtained by applying the four different formulations, share the same disturbance and initial state realization.

First, we focus on the share of constraint violations obtained with the four formulations.
This metric is computed by considering each point of each simulation independently.
% , since the controls are applied
% according to the receding horizon principle and the constraint \eqref{eq: explicit version of the chance constraint}
% is meant to be single chance in the respective formulations.
Therefore, the share of constraint violations is computed as:
\begin{equation}
    \frac{1}{N_\mathrm{s} \cdot T}\left[\sum_{k=0}^{N_\mathrm{s}-1} \sum_{t=0}^{T-1} \mathcal{I}_{> 0}(y_{k, t} - y_{\max})\right],
\end{equation}
where $\mathcal{I}_{> 0}$ denotes the indicator function that takes value $1$ if its argument is positive and $0$ otherwise.
Results are collected in Table~\ref{tab: monte carlo simulations - mean of violation and cost}.
The \textit{nominal} formulation is not aware about uncertainty and places the mean estimate $\hat{y}$ always exactly on the constraint.
Thus it is violated in about $50\%$ of the cases.
Instead, for the other formulations, violations happen on average $2.2-2.4\%$ of the time, corresponding to the chosen allowed constraint violation probability of $\epsilon \simeq 2.55\%$.
% Fig.~\ref{fig: monte carlo simulations - constraint violations} pictures the empirical probability density functions (PDF) based on the realizations of the uncorrupted output.

\begin{figure}
    \begin{center}
        \includegraphics[width=\linewidth]{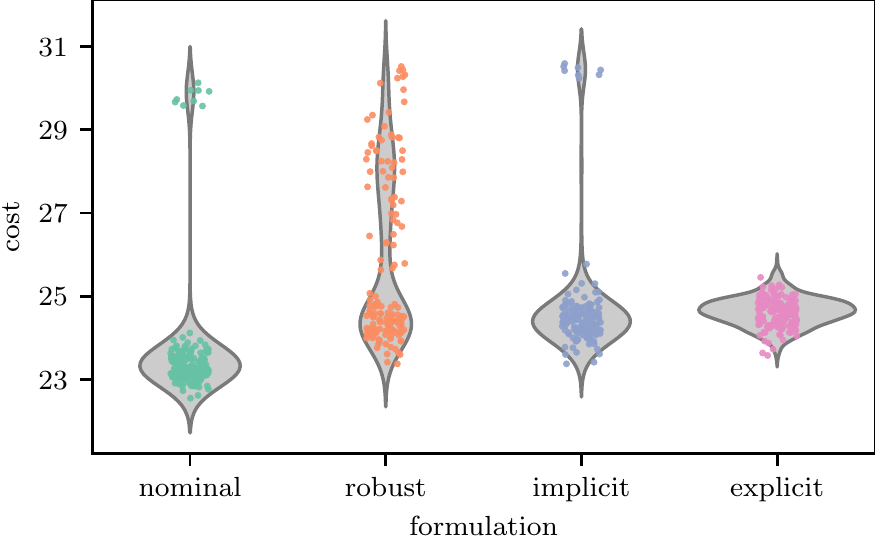}
    \end{center}
    \vspace*{-3mm}
    \caption{Extensive comparison -- Cost of the closed-loop trajectories \\
             The colored dots correspond to a portion of the total realizations and their scattering along the x-axis is to improve visualization.
             The shaded gray areas are the empirical probability density functions computed using all the available realizations.}\label{fig: monte carlo simulations - cost}
\end{figure}

Secondly, we take into account the cost of the closed-loop trajectories obtained with the different formulations.
In Table~\ref{tab: monte carlo simulations - mean of violation and cost}, we report the empirical mean of the cost distribution. One can see that the \textit{nominal} formulation leads to minimum cost, the \textit{robust} formulation to the highest cost and the two \textit{dual} formulations have similar cost which are in between the other two.
In Fig.~\ref{fig: monte carlo simulations - cost} we plot for each formulation the empirical PDF overlapped by a portion of the total realizations, specifically $15\%$ of the points.
The \textit{nominal} and \textit{implicit dual} formulations share a similar PDF where a cloud of points is far from the mean.
These realizations correspond to the scenarios where the scrap selection $u$ cannot drive the state estimate $\hat{x}$ close to the true state $x$, leading to bad performance.
However, these scenarios are few compared to the ones realized using the \textit{robust} formulation.
The latter has very scattered realizations leading to unreliable performance in terms of cost.
Finally, the \textit{explicit dual} formulation is outperforming the others in terms of reliability.
In fact the cost PDF does not have long tails resulting in very consistent behavior.

\begin{table}
    \vspace*{3mm}
    % \caption{Extensive comparison \\ Performance of the closed-loop trajectories}
    \caption{Performance of the closed-loop trajectories}
    \centering
    \ra{1.2}
    \begin{tabular}{@{}lrrrr@{}} \toprule
        Empirical mean & \textbf{\textit{Nominal}}& \textbf{\textit{Robust}}& \textbf{\textit{Impl. dual}}& \textbf{\textit{Expl. dual}} \\ \midrule
        Constraint viol. & 50.92\% & 2.43\% & 2.34\% & 2.18\% \\
        Cost           & 23.76 & 25.84 & 24.80 & 24.63 \\ \bottomrule
    \end{tabular}\label{tab: monte carlo simulations - mean of violation and cost}
\end{table}

% \begin{figure}
%     \begin{center}
%         \includegraphics[width=\linewidth]{img/mc_constraint_violations_cdf.pdf}
%     \end{center}
%     \vspace*{-3mm}
%     \caption{Extensive comparison -- Constraint realizations \\
%              The colored dots correspond to a portion of the total realizations and their scattering along the x-axis is to improve visualization.
%              The shaded gray areas are the empirical probability density functions computed using all the realizations.
%              To improve visualization, the y-axis is zoomed in.
%              From each formulation we left out $2.7\%$, $0.1\%$, $0.1\%$, $0.5\%$ of the points that violate the constraint respectively.}\label{fig: monte carlo simulations - constraint violations}
% \end{figure}
Finally, we want to motivate why the hyperparameter in the \textit{explicit dual} formulation is set to $\alpha=100$.
We compare five different values of $\alpha \in \{0, 1, 10, 100, 1000\}$.
The share of constraint violation and the corresponding distributions are very similar and therefore we focus on the cost.
The results are collected in Table~\ref{tab: monte carlo simulations - tuning alpha}.
We can see that the mean of the empirical PDF obtained with $\alpha=100$ achieves the minimum value and allows a more compact distribution of the cost with $99\%$ of the points below 25.53 as shown in the second column.
Also note that once the weight is set sufficiently large, i.e., $\alpha=10$, the results are very insensitive to the choice of $\alpha$.

In this work we have focused only on copper, but it is straightforward to extend the approach to multiple residual elements.

\begin{table}
    \caption{Comparing different $\alpha$ in the explicit dual formulation}
    \centering
    \ra{1.2}
    % \begin{tabular}{@{}lrr@{}} \toprule
    %     & \multicolumn{2}{c@{}}{\textbf{Cost}} \\ \cmidrule{2-3}
    %     $\alpha$ & Mean & Quantile 0.99 \\ \midrule
    %     0 & 24.80 & 30.60 \\
    %     1 & 25.07 & 30.80 \\
    %     10 & 24.85 & 26.23 \\
    %     100 & 24.63 & 25.53 \\
    %     1000 & 25.31 & 26.11 \\ \bottomrule
    % \end{tabular}
    % Preliminary version of the transposed table
    \begin{tabular}{@{}llrrrr@{}} \toprule
            \multicolumn{2}{r}{$\alpha$} & 1 & 10 & 100 & 1000 \\ \midrule
            \multirow{2}*{Cost} & Mean & 25.07 & 24.85 & 24.94 & 25.31\\
             & Quantile 0.99 & 30.80 & 26.23 & 25.80 & 26.11 \\ \bottomrule
        \end{tabular}\label{tab: monte carlo simulations - tuning alpha}
\end{table}

\section{Conclusions}
In this paper we have presented a formulation for both implicit and explicit dual model predictive control for a steel recycling process.
The process model allows for exact predictions of the state estimate uncertainty with a Kalman filter, leading to an uncertainty-aware and dual formulation.
The implicit formulation indirectly tackles the explore-exploit trade-off while in the explicit formulation one should balance the trade-off by tuning one hyperparameter.
The numerical simulations of the steel recycling process show that dual formulations outperform the uncertainty aware robust formulation in terms of cost of the closed-loop trajectories while achieving the same prescribed probability of constraint satisfaction.
Specifically, the explicit dual formulation provides a more consistent closed-loop behavior than the implicit dual formulation, and proves to be insensitive to the choice of the hyperparameter once it is set large enough.

\bibliography{syscop}

\end{document}